Original articles

# Different PCA approaches for vector functional time series with applications to resistive switching processes

C. Acal [a],[*], A.M. Aguilera [a], F.J. Alonso [a], J.E. Ruiz-Castro [a], J.B. Roldán [b]

[a] *Dept. of Statistics and O.R. and IMAG, University of Granada, Granada, Spain*
[b] *Dept. of Electronics and Computer Technology, University of Granada, Granada, Spain*

A B S T R A C T

This paper is motivated by modeling the cycle-to-cycle variability associated with the resistive switching operation behind memristors. Although the data generated by this stochastic process are by nature current–voltage curves associated with the creation (set process) and destruction (reset process) of a conductive filament, the statistical analysis is usually based on analyzing only the scalar time series related to the reset and set voltages/currents in consecutive cycles. As the data are by nature curves, functional principal component analysis is a suitable candidate to explain the main modes of variability associated with these processes. Taking into account this data-driven motivation, in this paper we propose two new forecasting approaches based on studying the sequential cross-dependence between and within a multivariate functional time series in terms of vector autoregressive modeling of the most explicative functional principal component scores. The main difference between the two methods lies in whether a univariate or multivariate PCA is performed so that we have a different set of principal component scores for each functional time series or the same one for all of them. Finally, the sample performance of the proposed methodologies is illustrated by an application on a bivariate functional time series of reset/set curves.

## 1. Introduction

The aim of this paper is to analyze functional time series (FTS) which consist of the observation of a series of dependent continuous functions instead of dependent real-valued data. From a mathematical viewpoint, we can define a scalar time series as $X_i$ with $i = 1, \ldots, n$ being the time points where the random variable is observed (also called index of the series). In contrast, a FTS is given by a set of dependent random functions $X_i(t)$ with $t$ in a real interval $T$. Functional time series can also be seen as a sequence of functional data originated by the partition of a continuous-time stochastic process in segments of equal length [21,22].

The classical Box-Jenkins methodology [11] is not adequate for modeling functional time series because the observations are random functions, whose observations are typically assumed to be elements of a separable Hilbert space [21]. Therefore, we are not dealing with vector spaces of finite dimension anymore. A solution is to consider some approach based on Functional Data Analysis (FDA), which comprehends multiple statistical methods for the analysis of functions (generally curves) that represent the evolution of a random variable over some continuous argument such as time or space [30,32]. The inherent problem in dealing with functional data is the infinite dimension of the space to which sample curves belong. For this reason, dimension reduction techniques such as Functional Principal Component Analysis (FPCA) are usually applied in FDA studies [3,9,14,29,42]. In particular, FPCA aims to

---

explain the main modes of variation of a functional dataset in terms of a small number of uncorrelated variables called principal components (see [38] for a comprehensive review).

The majority of FDA methodologies assume that the observed curves are generated by independent and identically distributed functional variables. This assumption is not true for functional data sequentially registered over time. For this reason, functional time series analysis have received special attention in last years. A deep review about the state-of-the-art of statistical time series models for functional data with sequential dependence can be consulted in [12]. Within the context of predicting FTS, the range of applications is very wide. A topic that is gaining interest recently is the dynamic updating. This consists of the forward prediction of a curve which is already being observed [39,40]. In relation to forecast several pieces of those curves that have not been totally observed, a new functional approach that uses all past observed curves is introduced in [25]. However, one of the most important issues in FTS analysis is the one-step-ahead forecasting, which has been widely addressed in the literature (see, e.g, [4,7,8,15,22,26,33]).

In the context of how FTS can be predicted by means of standard multivariate techniques, a simple and easily implantable algorithm is introduced in [8]. The underlying idea is to consider some approximation of the functional variable in a finite-dimension space generated by an orthonormal basis, such as for example the eigenfunctions of the covariance kernel, and to apply some multivariate forecasting methods over the vector-valued time series given by the basis coefficients. The objective of the present paper is to extend this idea to the general situation of a multivariate functional time series by analyzing the between and within linear dependence of several functional variables sequentially observed at the same time. As starting point we adapt the univariate functional autoregressive (FAR) model, deeply studied in [10] as the autoregressive Hilbertian process, to the multivariate case (vector functional autoregressive model, VFAR). Then, replacing the functional variables by their principal component representations (univariate or multivariate), the VFAR model is turned into a vector autoregressive model for the most explicative principal component scores of all considered functional time series.

This work is motivated by the analysis of data measured on Resistive Random-Access Memories (RRAMs). These electronic devices are compatible with the actual processes in the microelectronics industry and academy for their size, simplicity, and outstanding operational features. In fact, they are being currently incorporated in commercial chips [27]. These devices change their resistive state by creating (set) and destroying (reset) conductive filaments that short-circuit two metal electrodes separated by a dielectric. Then, the resistance state is fixed if the outside power is switched off. That is, the devices can easily switch and maintain two different resistance states, as a result, they are suitable for non-volatile memory applications [27].

The nature of the physical processes behind this resistive switching operation is stochastic; consequently, the measurements of their current in the laboratory present variability [1]. The statistical analysis of resistive memory measurements is essential to develop their different industrial applications (hardware cryptography [43], neuromorphic computing [44], etc.). There are different simulation and modeling schemes to describe the operation of RRAMs. One approach is based on detailed numerical simulations; for instance, kinetic Monte Carlo [5] and circuit-breaker-based tools [28]. RRAMs can also be described from a compact modeling perspective, linked to circuit design [13]. Another approach is connected to advanced statistical modeling of device measurements. In this sense, different statistical approaches based on phase-type distributions or functional data analysis have been employed (see, e.g. [2,36,37] and the references therein). These tools provide a coherent picture of RRAM operation due to their statistical description of the device inherent variability [35].

In this work, we have used current versus voltage reset and set curves (generated consecutively cycle-to-cycle within the resistive switching processes) because of their intrinsic variability and the interest in the context of industrial applications. An advanced perspective of functional data analysis is the viewpoint employed here. The main objective of this paper is to study the dependency relationships between set and reset processes to better understand the data structure and the device characteristics, which is essential for electronic circuit design. So far, the usual statistical analysis is based on modeling only the scalar time series related to the reset and set points (voltages/currents) in consecutive cycles [6,34]. This gives rise to misleading results since the functional nature of data is not taken into account during the analysis. In other words, important information about the device variability might be lost by considering scalar time series [35]. The FPCA approaches for the estimation of VFAR models proposed in this paper will help to improve the linear modeling and prediction of the bidimensional functional series of reset-set curves.

The rest of manuscript is organized as follows: the theoretical aspects related to the proposed methodology are in Section 2. Section 3 contains the application with data from resistive random access memories. The main conclusions are in Section 4.

## 2. Analysis of vector functional time series

A functional time series is a collection of dependent random functions denoted as $\{X_i(t) : i \in \mathbb{N}; t \in T = [a,b]\}$ where $i$ is the index of the series and $t$ is in some continuous domain $T$. This definition can be generalized for a vector of more than one functional variable so that $\mathbf{X}_i(t) = (X_{i1}(t), X_{i2}(t), \ldots, X_{iH}(t))^T$. Let us suppose that these vectors of functions are generated from a $H$-dimensional stochastic process whose components are second order and continuous in quadratic mean, with sample curves belonging to the Hilbert space $L^2[T]$ of squared integrable functions with the usual inner product.

### 2.1. Multivariate functional autoregressive model

The analysis of multivariate functional time series tries to study the sequential cross-dependence among several sets of curves with the aim of modeling their evolution over time and making predictions with enough precision.

The functional autoregressive model (FAR) under Hilbert space is likely the most popular technique for modeling an univariate FTS (see [10] for a detailed study). In this paper, the univariate FAR model of order $p$ is adapted for the multivariate case as follows (VFAR($p$)):

$$\mathbf{X}_{ij}(t) - \boldsymbol{\mu}_j(t) = \sum_{k=1}^{p}\sum_{h=1}^{H} \rho_{jh,k}(\mathbf{X}_{i-k,h}(t) - \boldsymbol{\mu}_h(t)) + \epsilon_{ij}(t) = \sum_{k=1}^{p}\sum_{h=1}^{H} \int_T \phi_{jh,k}(t-s)(\mathbf{X}_{i-k,h}(s) - \boldsymbol{\mu}_h(s))ds + \epsilon_{ij}(t), \quad (1)$$

where $\boldsymbol{\mu}(t) = (\mu_1(t), \ldots, \mu_H(t))^T$ is the mean function vector of the FTS $\mathbf{X}_i(t)$, $\phi_{jh,k}(t)$ is the kernel function of the bounded linear Hilbert–Schmidt operator $\rho_{jh,k}$ from $L^2[T]$ to $L^2[T]$, $\mathbf{X}_{i-k}$ represents the $k$th lag of $\mathbf{X}_i$ and $\epsilon_i(t)$ denotes a $L^2[T]$-vector white noise with zero mean and finite second moments. The above model assume stationarity which means the dynamics of each series is stable over time. A revision of models under nonstationarity hypothesis can be seen in [12].

Let us denote by $\rho_k = (\rho_{jh,k})$ the matrix of dimension $H \times H$ whose elements are the operators representing the serial cross-dependence between the different pairs of functional variables on their past values. Then, the VFAR(p) model in Eq. (1) can be expressed in matrix form as follows:

$$\mathbf{X}_i(t) - \boldsymbol{\mu}(t) = \sum_{k=1}^{p} \rho_k(\mathbf{X}_{i-k}(t) - \boldsymbol{\mu}(t)) + \boldsymbol{\epsilon}_i(t) = \sum_{k=1}^{p} \int_T \phi_k(t-s)(\mathbf{X}_{i-k}(s) - \boldsymbol{\mu}(s))ds + \boldsymbol{\epsilon}_i(t), \quad (2)$$

with $\phi_k = (\phi_{jh,k})$ being the $H \times H$ matrix of kernel functions associated with the operators in $\rho_k$.

The classical estimation methods such us maximum likelihood do not work properly with FAR models due to the infinite dimension of the data. In this paper, we propose two estimation approaches for VFAR models based on the usual projection of functional data on a finite-dimension orthogonal basis of functions. Specifically, we will consider the basis of eigenfunctions of the covariance and cross-covariance operators of the processes generating the FTS and their associated FPCA decompositions.

### 2.2. Functional principal component analysis

Let $\mathbf{X}_1(t), \ldots, \mathbf{X}_n(t)$ be a sample of vector functions where $\mathbf{X}_i(t) = (X_{i1}(t), \ldots, X_{iH}(t))^T$ contains the information for the $i$th subject on each functional variable $X_h$, $h = 1, \ldots, H$. Besides, let us assume that this sample of vector curves are realizations of a multivariate stochastic process with sample vector mean function $\boldsymbol{\mu} = (\mu_1(t), \ldots, \mu_H(t))^T$ and sample matrix covariance function $\mathbf{C}$ such that $\mathbf{C}(t,s) = (C_{h,h'}(t,s))$, $t, s \in T = [a,b]; h, h' = 1, \ldots, H$

$$C_{h,h'}(t,s) = \frac{1}{n}\sum_{i=1}^{n} \left(X_{ih}(t) - \mu_h(t)\right)\left(X_{ih'}(s) - \mu'_h(s)\right).$$

Let us observed that $C_{h,h'}$ will be the covariance function denoted by $C_h$ if $h = h'$ and a cross-covariance function otherwise.

#### 2.2.1. Univariate FPCA

It is well known that Functional PCA aims to explain the main modes of variation of a sample of curves in terms of a reduced number of centered uncorrelated scalar variables with maximum variance called functional principal component scores (PCs) [31,32].

For each of the considered functional variables $X_h(h = 1, \ldots, H)$, the $j$th principal component is defined as a generalized linear combination of the functional variable

$$\xi_{ij}^h = \int_T X_{ih}(t) f_j^h(t)\, dt, \ i = 1, \ldots, n,$$

with the weight function $f_j^h$ obtained by maximizing the variance $Var\left[\int_T X_{ih}(t) f(t) dt\right]$ r.t. $\|f\|^2 = 1; \int f_k^h(t) f(t) dt = 0, k = 1, \ldots, j-1$.

The solutions to this optimization problem are given by the eigenfunctions of the sample covariance operator, $C_h f_j^h = \lambda_j^h f_j^h$, with $\{\lambda_j^h\}_{j \geq 1}$ being the decreasing sequence of non-null eigenvalues such that $\lambda_j^h = Var[\xi_j^h]$ and $C_h(f(s)) = \int C_h(t,s) f(t) dt$.

Then, the well-known Karhunem–Loéve (K-L) expansion provides the following approximated principal component decomposition for each functional variable $X_h$ in terms of the first $q_h$ principal components for which the proportion of cumulative variability is as close as possible to one [14]:

$$X_{ih}^{q_h}(t) = \mu_h(t) + \sum_{j=1}^{q_h} \xi_{ij}^h f_j^h(t) = \mu_h(t) + (\boldsymbol{F^h}(t))^T \boldsymbol{\xi_i^h}, \quad (3)$$

with $\boldsymbol{F^h}(t) = (f_1^h(t), \ldots, f_{q_h}^h(t))^T$ and $\boldsymbol{\xi_i^h} = (\xi_{i1}^h, \ldots, \xi_{iq_h}^h)^T$ being the vectors of the first $q_h$ principal components weight functions and scores, respectively.

#### 2.2.2. Multivariate FPCA

Multivariate Functional Principal Component Analysis (MFPCA) is the natural extension of the univariate FPCA when there are more than one functional variable in the analysis.

The $j$th principal component is computed by means of the following expression (see [24] for a detailed explanation):

$$\xi_{ij} = \int_T (\boldsymbol{X_i}(t) - \boldsymbol{\mu}(t))^T \boldsymbol{f_j}(t)dt = \sum_{h=1}^{H} \int_T (X_{ih}(t) - \mu_h(t))f_{jh}(t)dt,$$

where the weight functions $\boldsymbol{f_j}(t) = (f_{j1}(t), \ldots, f_{jH}(t))^T$ are the solutions to the eigenequation system $C\boldsymbol{f_j} = \lambda_j \boldsymbol{f_j}$ that can be expressed as

$$\sum_{h'=1}^{H} \int_T C_{hh'}(t,s) f_{jh'}(t) dt = \lambda_j f_{jh}(s), \forall h = 1, \ldots, H.$$

Note that $\lambda_j = \text{Var}[\xi_j]$ and finally, the multivariate functional variable can be approximated by truncating the K-L expansion in terms of the first $q$ principal components as

$$\boldsymbol{X}_i^q(t) = \boldsymbol{\mu}(t) + \sum_{j=1}^{q} \xi_{ij} \boldsymbol{f_j}(t) = \boldsymbol{\mu}(t) + \boldsymbol{F}(t)\boldsymbol{\xi_i}, \quad (4)$$

with $F(t)$ being a $H \times q$ matrix whose columns form the orthonormal basis of vector eigenfunctions $\boldsymbol{f_j}(t)$ and $\boldsymbol{\xi_i} = (\xi_{i1}, \ldots, \xi_{iq})^T$ the vector of the first $q$ principal component scores.

As in the case of univariate FPCA, in practice $q$ must be chosen so that the percentage of explained cumulative variance is high enough (around 95% at least) in order to guarantee a good reconstruction of the sample curves.

### 2.3. New functional principal component vector autoregressive models

Reducing the infinite dimension of the vector functional time series by using FPCA (univariate or multivariate) and replacing the functional variables by their univariate or multivariate principal component representation (Eqs. (3) or (4)), the VFAR(p) model expressed by Eqs. (1) and (2) is equivalent to a classical vector autoregressive model (VAR(p)) for the vector time series with the most explicative principal component scores

$$\boldsymbol{\xi_i} = \sum_{k=1}^{p} \boldsymbol{\Omega_k} \boldsymbol{\xi_{i-k}} + \boldsymbol{\epsilon_i^*},$$

where $\boldsymbol{\epsilon_i^*}$ is a vector white noise, $\boldsymbol{\Omega_k}$ is the matrix of coefficients of order $q \times q$ and $\boldsymbol{\xi_i}$ is a principal component vector time series of dimension $q$ defined in different ways according to the univariate or multivariate case as follows:

- Univariate functional PC: $\boldsymbol{\xi_i} = (\xi_1^1, \ldots, \xi_{q_1}^1, \ldots \xi_1^h, \ldots, \xi_{q_h}^h, \ldots \xi_1^H, \ldots, \xi_{q_H}^H)^T$, with $q = \sum_{h=1}^{H} q_h$. That means that the first $q_h$ PC scores are considered for each functional time series $X_{ih}, h = 1, \ldots, H$.
- Multivariate functional PC: $\boldsymbol{\xi_i} = (\xi_1, \ldots, \xi_q)^T$. That means that the common first $q$ PC scores of the multivariate functional time series are selected to fit the VAR(p) model.

Let us observe that the initial problem would be reduced to estimate one of these vector autoregressive models for the most explicative principal components. Then, predictions for the multivariate FTS are given by the following functional principal component autoregressive models of order $p$:

- Multivariate FPCA approach: MFPCA-VAR(p)

$$\boldsymbol{X}_i^q(t) = \boldsymbol{\mu}(t) + \sum_{k=1}^{p} \boldsymbol{F}(t)\boldsymbol{\Omega_k}\boldsymbol{\xi_{i-k}} + \boldsymbol{F}(t)\boldsymbol{\epsilon_i^*}. \quad (5)$$

Let us observe that the matrix of operators $\rho_k$ in Eq. (2) is estimated by the linear transformation, through coefficient matrix $\boldsymbol{\Omega_k}$, of the coordinates on the space spanned by the orthonormal basis $F$. That is, $\forall g \in (L^2[T])^H$ so that $g(t) = \boldsymbol{F}(t)\boldsymbol{g}$ then $\rho_k(g(t)) = \boldsymbol{F}(t)\boldsymbol{\Omega_k}\boldsymbol{g}$.

- Univariate FPCA approach: FPCA-VAR(p)
  In this case we have a prediction model for each functional time series

$$X_{ih}^{q_h}(t) = \mu_h(t) + \sum_{k=1}^{p} (\boldsymbol{F^h}(t))^T \boldsymbol{\Omega_k^h} \boldsymbol{\xi_{i-k}}, \quad (6)$$

where $\boldsymbol{\Omega_k^h}$ is the $q_h \times q$ matrix of coefficients associated with the PCs of the $h$th functional time series: $\xi_i^h = \sum_{k=1}^{p} \boldsymbol{\Omega_k^h} \boldsymbol{\xi_{i-k}}$.

Summarizing, The VFAR methodology proposed in this paper is based in a two step FDA approach:

1. Estimating FPCA of the two involved functional variables (reset and set processes). Two different methods are considered depending on whether univariate or multivariate FPCA is performed.
2. VAR modeling of the most explicative functional principal component scores.

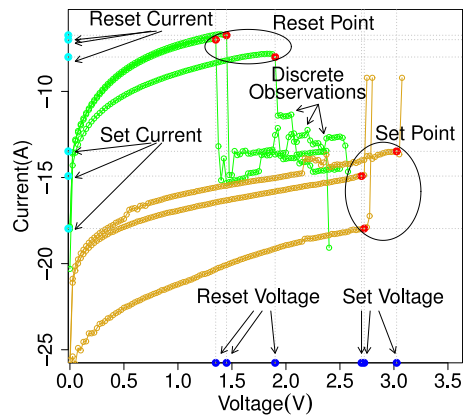

**Fig. 1.** Napierian logarithm of the experimental current versus voltage for several curves of continuous resistive switching cycles (set and reset curves in orange and green, respectively). Set/reset points are indicated for clearness. Discrete observations are shown with circles in each path. (For interpretation of the references to color in this figure legend, the reader is referred to the web version of this article.)

The asymptotic and finite-sample properties for the estimators of the statistical models considered in our approach have been deeply studied in the literature. Asymptotic theory for FPCA of a Hilbert random variable and applications in statistical inference were deeply studied in [14]. The finite-sample performance of the approximate estimation of FPCA from B-spline smoothing of the sample curves with and without penalization was studied by simulation in [3]. On the other hand, the asymptotic properties and sample behavior of multivariate FPCA has been recently investigated in [20]. With respect to the VAR models considered for modeling the most explicative PC scores, the theoretical properties of coefficient estimators have been thoroughly analyzed in the literature. It has been proven that the estimates of the coefficient matrices in a VAR model (estimation performed using OLS formulas) are consistent even when the innovations are not normally distributed (see Proposition 11.1 in [19]). In other words, the consistency of the coefficient estimates in a VAR model is a large-sample property that does not depend on the distribution of the innovations, as long as certain regularity conditions are met. On the other hand, the univariate scenario in which only one functional time series is modeled and predicted, was already investigated both theoretically and via simulation by [8], whose starting point was the univariate functional autoregressive model, also deeply studied by [10]. This methodology has been extended in the current paper for the case of multivariate functional time series with the aim of modeling the internal relationships of set and reset processes.

## 3. Application

The proposed FTS approaches are applied to model and forecast the bivariate functional time series of the reset/set curves collected from the resistive switching processes behind RRAMs devices. These models can help to a better design and simulation of these devices which are currently great challenges in the area of microelectronic.

### 3.1. Data description

The device measurements are obtained from RRAMs with nickel and silicon electrodes and a dielectric 10 nm thick of hafnium oxide. The fabrication details were given in [17]. A semiconductor parameter analyzer was employed to obtain a long series of resistive switching cycles by means of external programmed triangular voltage signals. The silicon electrode was grounded, and a negative voltage was applied to the nickel electrode, although henceforth we assume the absolute value of applied voltage for the sake of simplicity (see Fig. 1). Several set (orange symbols) and reset (green symbols) current versus voltage curves of a series of measured resistive switching cycles (1826 cycles) are plotted in Fig. 1. A sudden current drop corresponds to the rupture of the conductive filament (known as reset point) that short-circuits the electrodes (see that the conductive filament works here as a fuse). The corresponding voltages and currents are known as reset voltages and currents respectively, which have been explicitly shown in Fig. 1. In set curves a sudden current rise takes place when the conductive filament is fully formed, reducing drastically the device conductance. At this point the set voltage and current are defined (see Fig. 1).

We have determined the set and reset points in all curves (corresponding to 1826 consecutive cycles) by considering a current drop above 20% in two consecutive voltage points. Hereinafter, the objective is to analyze the curves until the corresponding set/reset point.

### 3.2. Data approximation by FPCA

Several problems are found when functional data analysis is applied for modeling experimental data measured on RRAM memories. On the one hand, there is not a common domain for curves because of set and reset points are different for each cycle.

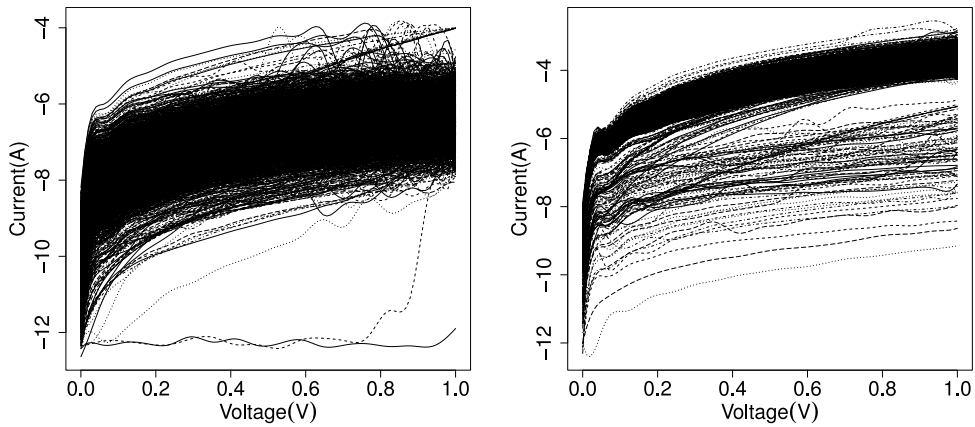

**Fig. 2.** B-spline smoothing of set (left panel) and reset (right panel) curves. The napierian logarithm of the experimental current is displayed on axis Y.

**Table 1**
Cumulative variability (%) explained by the first six principal components.

| PC    | 1       | 2       | 3       | 4       | 5       | 6       |
|-------|---------|---------|---------|---------|---------|---------|
| Reset | 91.6264 | 96.7037 | 98.9527 | 99.4620 | 99.8493 | 99.9275 |
| Set   | 83.6191 | 90.7955 | 93.9711 | 95.6772 | 96.7411 | 97.4494 |
| Mult. | 75.4757 | 85.7073 | 91.1773 | 93.7822 | 95.2890 | 96.2322 |

On the other hand, given the inability to compile the full information of curves in practice, multiple discrete observations at a finite set of current values until the set/reset point are recorder for each curve (see Fig. 1).

To solve these issues, an option is to follow the procedure described in [2]. First, all curves are registered in the interval [0, 1] by dividing the discrete observations by the corresponding set/reset point. Second, the functional form of curves is reconstructed by assuming that they belong to a finite-dimension space generated by a basis [31,32]. Several procedures have been considered in the literature for selecting the optimum dimension of the basis such as cross-validation, generalized cross-validation with a roughness penalty or more recently, Bayesian approaches via Gibbs sampler [16,32,41]. Following an ad-hoc graphical inspection to avoid the typical overfittings reached by means of numerical procedures, a cubic B-spline basis of dimension 20 is considered in this paper for smoothing the sample curves. Fig. 2 reveals that several curves could be outliers in shape and/or in range. This fact must be studied because FPCA is highly sensitive in the presence of outliers. To identify possible anomalous curves, the functional bagplot method proposed by Hyndman and Shang, which is based on Tukey's halfspace location depths, is employed (see [23] for more detail). After this analysis, 7% of cycles were removed. This decision makes sense from a physical viewpoint because a curve identified as outlier means that the conductive filament has not been totally formed (or destructed for the case of reset curves). That is, there has been a measurement error and therefore, the entire cycle must be eliminated.

Next, FPCA and MFPCA are applied to the smoothed registered curves. Hereinafter, last ten cycles are excluded from the analysis to constitute a test sample that will be used to make predictions and check the quality of the fitted models. The percentage of cumulative variability explained by the first six principal component scores for both univariate and multivariate case can be seen in Table 1.

### 3.3. Time series modeling of PC

The VAR modeling of functional principal components for set and reset trajectories can be done following two alternative ways. The first consists of time series analysis on the PCs obtained by multivariate FPCA and the second is based on the univariate FPCA of set and reset variables separately. In both approaches, the selected components must exceed 95% of total variability in order to guarantee a good reconstruction of the underlying processes. Thus, the first five PCs are taken for fitting the MFPCA-VAR model and the first two and four PCs for the reset and set processes, respectively, in the univariate setting.

The first step of the TS analysis lies in establishing the causal relationships among components. For this purpose, a suitable option is the Granger's causality study by comparing two-to-two the selected PCs [18]. Granger's test consists of analyzing the significance of the coefficients by comparing nested models. Selected two components $\xi_{j'}$ and $\xi_j$, the following model is defined to study if $\xi_j$ causes $\xi_{j'}$:

$$\xi_{ij'} = \sum_{l=1}^{p} \phi_l \xi_{i-l,j'} + \sum_{k=1}^{q} \phi'_k \xi_{i-k,j},$$

so that the analysis is reduced to test if the parameters $\phi'_k$ are significantly different from zero. Rejecting the null-hypothesis ($H_0 : \phi'_k = 0, k = 1, \ldots, q$) would imply that $\xi_j$ causes $\xi_{j'}$. In broad terms, this test checks if the variable $\xi_j$ contains relevant information for modeling the variable $\xi_{j'}$. Otherwise, the past of $\xi_{j'}$ would be sufficient. The inherent problem is to determine the values $p$ and $q$ in the model. This issue can be solved by fitting univariate models for each component and studying the Granger's causality over the residuals.

Table 2
Granger-causality test over residuals of univariate models for multivariate PCs.

| | MPC1 | MPC2 | MPC3 | MPC4 | MPC5 |
|---|---|---|---|---|---|
| MPC1 | | ← | ← | | ← |
| MPC2 | | | ← | | |
| MPC3 | ← | ← | | ← | ← |
| MPC4 | | ← | ← | | |
| MPC5 | | | ← | | |

Table 3
Granger-causality test over residuals of univariate models for univariate PCs.

| | RPC1 | RPC2 | SPC1 | SPC2 | SPC3 | SPC4 |
|---|---|---|---|---|---|---|
| RPC1 | | ← | | | | |
| RPC2 | ← | | ← | | | ← |
| SPC1 | ← | ← | | ← | | ← |
| SPC2 | | ← | ← | | | ← |
| SPC3 | ← | | | | | |
| SPC4 | | ← | | ← | | |

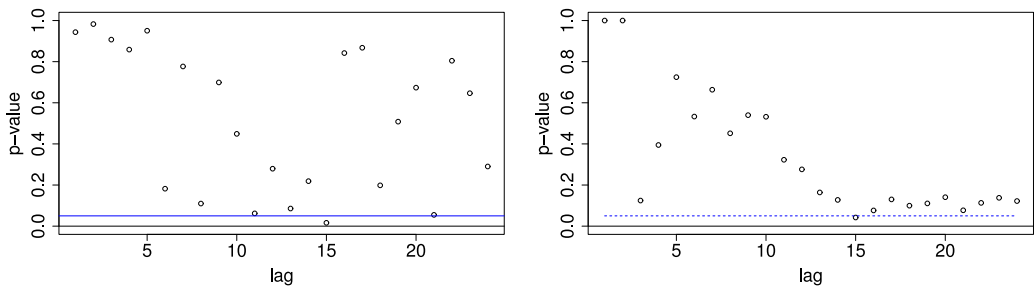

**Fig. 3.** P-values for the study of the significance of the cross-correlation matrices of the residuals (left panel) and the Ljung–Box tests on the residuals (right panel).

Let us begin by the multivariate FPCA approach. The following ARMA models have been fitted for the first five principal components: ARMA(1,1), ARMA(1,1), ARMA(1,2), ARMA(1,2) and ARMA(0,2). The residuals of each one of these univariate models are regressed on the other residuals to test whether their coefficients are significantly different from zero. If they are significant, it will be said that there is causality of the independent variable over the dependent variable. A summary of the causality study over the residuals can be seen in Table 2, where, for example, the arrow in the last row indicates that the third component causes the fifth component.

Facing these results, a joint vector autoregressive modeling of all PCs must be carried out, since the dependence cannot be studied on vector series of lower dimension. Let us observed that the third PC causes all other PCs and PC3 is caused by all the others. Then, a VAR model is fitted for the vector time series of dimension five with all the PCs. The model order estimated by selecting the model with less AIC is seven. Non-significant coefficients are eliminated by setting the threshold 1.96 in the critical value used for studying the significance of coefficients.

We also performed an uncorrelatedness study of the five dimension model errors $\varepsilon_t$. If we study individually the significance of cross-correlation matrices (CCMs) between errors corr$(\varepsilon_t, \varepsilon_{t+k})$, it is desirable that these correlation matrices are null for values $k > 0$. However, the same as with the autocorrelation function in univariate series, some values may be significant due to randomness. For the model to fit well, it is necessary that there are not several significant CCMs and the first k are not significant ($k = 5$). In this case, only matrix 15 in the CCM plot shown in Fig. 3 is below the 0.05 band. Besides, if we look at the Ljung–Box plot which tests the significance of the first k CCMs jointly (the null hypothesis is that they are all null), only the $k = 15$ case is significant (see Fig. 3). Therefore, the fitted model can be considered as adequate.

The estimated VAR(7) model can be written as

$$MPC1_t = 0.348\,MPC1_{t-1} - 0.775\,MPC2_{t-1} + 0.169\,MPC1_{t-2} - 0.214\,MPC3_{t-2} +$$
$$+ 0.067\,MPC1_{t-3} + 0.074\,MPC1_{t-4} + 0.117\,MPC2_{t-5} - 0.170\,MPC3_{t-5} +$$
$$+ 0.066\,MPC1_{t-6} + \varepsilon_{1t}$$
$$MPC2_t = 0.134\,MPC2_{t-1} + 0.148\,MPC2_{t-2} - 0.093\,MPC3_{t-2} + 0.097\,MPC2_{t-3} -$$

**Table 4**
Final causality relationships among PC residuals after a partial causality study.

|      | RPC1 | RPC2 | SPC1 | SPC2 | SPC3 | SPC4 |
|------|------|------|------|------|------|------|
| RPC1 |      | ←    |      |      |      |      |
| RPC2 | ←    |      |      |      |      |      |
| SPC1 | ←    | ←    |      | ←    |      | ←    |
| SPC2 |      |      | ←    |      |      | ←    |
| SPC3 | ←    |      |      |      |      |      |
| SPC4 |      |      |      | ←    |      |      |

$$-0.030\,MPC1_{t-4} + 0.073\,MPC2_{t-4} + 0.152\,MPC5_{t-4} + 0.097\,MPC2_{t-6} + \varepsilon_{2t}$$

$$\begin{aligned}
MPC3_t =\ & -0.052\,MPC2_{t-1} + 0.208\,MPC3_{t-1} + 0.092\,MPC4_{t-1} + 0.111\,MPC5_{t-1} + \\
& + 0.090\,MPC3_{t-2} + 0.087\,MPC3_{t-3} - 0.021\,MPC1_{t-4} + 0.064\,MPC3_{t-4} + \\
& + 0.083\,MPC5_{t-4} + 0.082\,MPC3_{t-6} + 0.097\,MPC3_{t-7} + \varepsilon_{3t}
\end{aligned}$$

$$\begin{aligned}
MPC4_t =\ & -0.038\,MPC2_{t-1} + 0.039\,MPC3_{t-1} + 0.142\,MPC4_{t-1} + 0.049\,MPC4_{t-2} - \\
& - 0.065\,MPC5_{t-2} + 0.076\,MPC3_{t-4} + 0.091\,MPC4_{t-4} + 0.055\,MPC4_{t-5} - \\
& - 0.081\,MPC5_{t-7} + \varepsilon_{4t}
\end{aligned}$$

$$MPC5_t = 0.033\,MPC3_{t-1} + 0.057\,MPC5_{t-1} + 0.074\,MPC5_{t-2} + 0.036\,MPC3_{t-3} + \varepsilon_{5t}$$

Let us now consider the univariate design with the first 2 PCs of reset curves and the first 4 PCs of set curves. As in the multivariate case, we performed a causality study on the residuals of the univariate TS models fitted to the dimension six series of PCs. In the two-to-two causality study we obtain the results shown in Table 3.

Here, we also examined what we call partial causality. Partial causality is defined, in a multivariate context, as the causality of one variable on another given a third or group of variables. For example, we can study the partial causality of SPC1 on RPC2 given RPC1. RPC2 is regressed on the past values of RPC1 and SPC1 and compared with the regression of RPC2 on the past values of RPC1. A Wald test on the coefficients of the past values of SPC1 will indicate whether there is partial causality. In this case, we conclude that SPC1 does not partially cause RPC2 given RPC1. After a general revision, the final causality study is summarized in Table 4. Note that a partial causality analysis was also carried out for the multivariate setting (joint reset-set modeling) using five components. However, the causality structure was not simplified for that case.

From the causality structure of Table 4, it is deduced that a multivariate modeling of RPC = (RPC1, RPC2) and SPC = (SPC1, SPC2, SPC4) can be performed. Besides, a transfer function model is estimated for the third component of set process SPC3. Finally, an improvement can also be proposed by relating the residuals of SPC1 and those of the multivariate reset modeling.

Thus, the proposed models after readjusting them by eliminating non-significant coefficients are

- A VAR(9) model for the reset components

$$\begin{aligned}
\text{RPC}_t =\ & \begin{pmatrix} 0.269 & 0.352 \\ 0.009 & 0.514 \end{pmatrix} \text{RPC}_{t-1} + \begin{pmatrix} 0.137 & 0.000 \\ 0.000 & 0.201 \end{pmatrix} \text{RPC}_{t-2} + \\
& + \begin{pmatrix} 0.078 & 0.000 \\ 0.000 & 0.082 \end{pmatrix} \text{RPC}_{t-3} + \begin{pmatrix} 0.074 & 0.000 \\ 0.000 & 0.063 \end{pmatrix} \text{RPC}_{t-4} + \\
& + \begin{pmatrix} 0.067 & -0.358 \\ 0.000 & 0.000 \end{pmatrix} \text{RPC}_{t-6} + \begin{pmatrix} 0.053 & 0 \\ 0.000 & 0 \end{pmatrix} \text{RPC}_{t-7} + \\
& + \begin{pmatrix} 0.082 & 0.000 \\ 0.000 & 0.062 \end{pmatrix} \text{RPC}_{t-9} + \begin{pmatrix} \varepsilon_{1,t} \\ \varepsilon_{2,t} \end{pmatrix}.
\end{aligned}$$

- A VAR(8) model for the set components

$$\begin{aligned}
\text{SPC}_t =\ & \begin{pmatrix} 0.329 & 0.669 & 0.000 \\ 0.020 & 0.145 & 0.124 \\ 0.000 & 0.033 & 0.073 \end{pmatrix} \text{SPC}_{t-1} + \begin{pmatrix} 0.164 & 0.000 & 0.000 \\ 0.000 & 0.066 & 0.000 \\ 0.000 & 0.000 & 0.068 \end{pmatrix} \text{SPC}_{t-2} \\
& + \begin{pmatrix} 0.060 & 0.000 & 0.000 \\ 0.000 & 0.085 & 0.000 \\ 0.000 & 0.039 & 0.000 \end{pmatrix} \text{SPC}_{t-3} + \begin{pmatrix} 0.061 & 0.000 & -0.296 \\ 0.025 & 0.051 & 0.161 \\ 0.000 & 0.000 & 0.000 \end{pmatrix} \text{SPC}_{t-4} \\
& + \begin{pmatrix} 0.059 & 0.000 & 0.000 \\ 0.0000 & 0.082 & 0.000 \\ 0.000 & 0.000 & 0.000 \end{pmatrix} \text{SPC}_{t-6} + \begin{pmatrix} 0.081 & 0.000 & 0.000 \\ 0.000 & 0.071 & 0.000 \\ 0.000 & 0.000 & 0.000 \end{pmatrix} \text{SPC}_{t-8} + \begin{pmatrix} e_t \\ \varepsilon_{4,t} \\ \varepsilon_{6,t} \end{pmatrix},
\end{aligned}$$

with a transfer function model with two inputs for the errors of SPC1

$$e_t = 0.001 - 0.151\varepsilon_{1,t} - 0.612\varepsilon_{1,t-1} - 1.169\varepsilon_{2,t} + 0.923\varepsilon_{2,t-1} + \varepsilon_{3,t}.$$

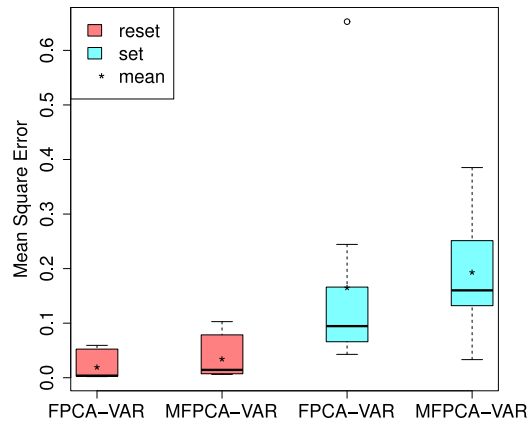

**Fig. 4.** Box-plots of integrated mean square error on the test sample.

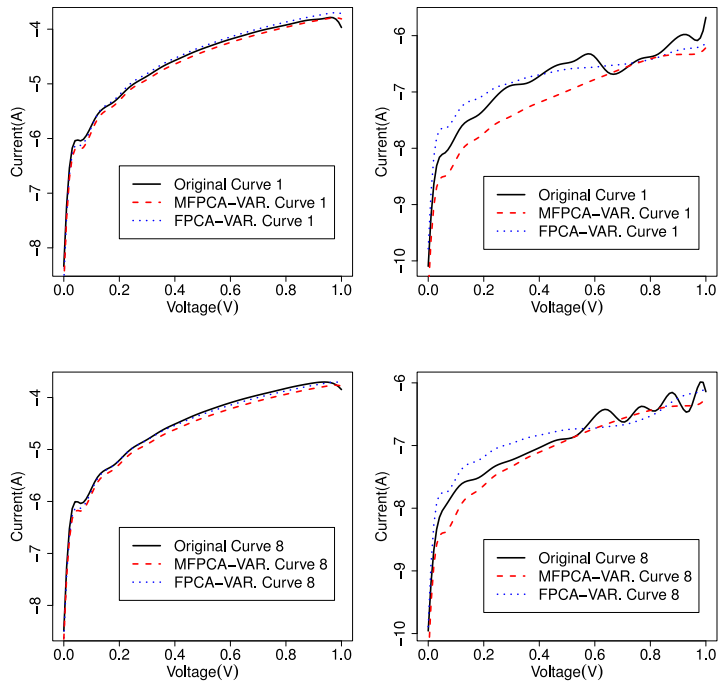

**Fig. 5.** Predictions for some cycles of the test sample: reset curves on the left side and set curves on right one. Original curves in solid-black line and the predictions with the univariate and multivariate approaches in dotted-blue line and dashed-red line, respectively. The napierian logarithm of the experimental current is displayed on axis Y. (For interpretation of the references to color in this figure legend, the reader is referred to the web version of this article.)

- A transfer function model for SPC3 with input RPC1

$$\text{SPC3}_t = 0.002 - 0.147\,\text{RPC1}_t - 0.440\,\text{RPC1}_{t-1} - 0.216\,\text{RPC1}_{t-2} + \frac{1 + 0.847 B + 0.057 B^2}{1 - 0.958}\varepsilon_{5,t},$$

with $B$ being the backward shift operator.

### 3.4. Predictions

Finally, we perform the prediction for the curves of the test sample. The vector models previously fitted are transformed to a functional time series prediction model using the truncated K-L expansions as in Eqs. (5) and (6). Theoretically, we are making predictions of the approximation of the original curves given by the PC decompositions which allow to describe the main characteristics of RRAMs when the cumulative percentage of explained variability is high enough (around 95% as in this paper). Therefore, we assume that the fitting by means of this orthogonal decomposition is adequate and our predictions should be looked like to the original cycles if the fitted models are good. The box-plots for the distribution of integrated mean square prediction error on the test sample curves appear in Fig. 4. On the other hand, predictions for some cycles of the test sample can be seen in Fig. 5.

Both the univariate and the multivariate approaches provide satisfactory results, especially for the reset process, with the univariate FPCA-VAR model having a little more complexity with a smaller dimension reduction and providing a bigger prediction accuracy. In the case of the set curves where the variability is higher (see Fig. 2), the forecasts do not control the peaks but are able to adequately estimate the trend. This fact suggests that increasing the number of PCs in the study would be a suitable solution to control all patterns of variation in the data. However, this noise, which is well-known and assumed due to the internal physical behavior of electrodes during the resistive switching operation, is usually ignored in the analysis for the experts in this area. Therefore, the new methodology is adequate to model the resistive switching stochastic processes and can be used to improve circuit design at the industrial level and to understand the characteristics and structure of the formation and rupture of the conductive filament.

## 4. Conclusions

The data-driven motivation of the current paper is to model the variability associated with the resistive switching processes behind the RRAMs operation. From the mathematical viewpoint, the data generated in these processes are a sample of bidimensional current–voltage curves corresponding to the consecutive cycles of formation (set) and rupture (reset) of a conductive filament. Given the stochastic and functional nature of data, we propose a new methodology based on functional data analysis for predicting new cycles and/or analyzing the sequential cross-dependence between and within set and reset processes through a causality study. In particular, we extend the univariate functional autoregressive model for the multivariate case and propose two estimation approaches based on vector autoregressive modeling for the vector-valued time series of the most explicative univariate or multivariate functional principal components. The application results show a better predictive behavior of the univariate FPCA-VAR model that is also able to detect causality relationships between the different functional series under study.

**CRediT authorship contribution statement**

**C. Acal:** Developed all the statistical theory, Computationally implemented the methodology, Obtained the results, Contributed equally to the writing part of the manuscript. **A.M. Aguilera:** Developed all the statistical theory, Computationally implemented the methodology, Obtained the results, Contributed equally to the writing part of the manuscript. **F.J. Alonso:** Developed all the statistical theory, Computationally implemented the methodology, Obtained the results, Contributed equally to the writing part of the manuscript. **J.E. Ruiz-Castro:** Developed all the statistical theory, Computationally implemented the methodology, Obtained the results, Contributed equally to the writing part of the manuscript. **J.B. Roldán:** In charge of the physical theoretical explanation and its conclusions, Performed the data curation for this study, Contributed equally to the writing part of the manuscript.

**Acknowledgments**

The authors acknowledge financial support by project PID2020-113961GBI00 of the Spanish Ministry of Science and Innovation (also supported by the FEDER programme) and projects PID2021-128077NB-I00 and PID2022-139586NB-44 funded by MCIN/AEI/10.13039/501100011033 and FEDER, UE. This research has been also supported by project FQM-307 of the Government of Andalusia (Spain) and the IMAG María de Maeztu grant CEX2020-001105-M/AEI/10.13039/501100011033. Funding for open access charge: Universidad de Granada / CBUA. All authors have read and agreed to the published version of the manuscript.